\documentclass[a4paper,12pt]{article}
\usepackage{amsmath,amsthm,amsfonts}
\usepackage{a4wide}
\usepackage[T1]{fontenc}
\usepackage{times}
\usepackage[hidelinks]{hyperref}

\DeclareMathOperator{\tr}{Tr}
\newcommand{\E}{\mathbb{E}}
\newcommand{\e}{\mathrm{e}}
\newcommand{\R}{\mathbb{R}}
\DeclareMathOperator{\cov}{Cov}
\DeclareMathOperator{\var}{Var}
\newcommand{\id}{\mathrm{Id}}
\renewcommand{\SS}{\mathbb{S}}
\renewcommand{\phi}{\varphi}

\newtheorem{theorem}{Theorem}

\theoremstyle{definition}
\newtheorem{remark}[theorem]{Remark}

\title{A note on slicing}
\author{Boaz Klartag and Joseph Lehec}
\date{}

\begin{document}

\maketitle

\begin{abstract}
We give a short proof of Bourgain's slicing theorem with the explicit
bound $L_K\leq\sqrt{5/\pi}$ for every convex body $K$. The argument
combines the strong thin-shell inequality obtained from the moment-measure
method of Chen and Klartag with an entropy estimate along Eldan's
stochastic localization. The proof uses neither $M$-ellipsoid theory
nor the reverse Santal\'o inequality.
\end{abstract}

\section{Introduction}

Bourgain's slicing problem asks whether there exists a universal constant
$c>0$ such that every convex body $K\subset\R^n$ of volume one has a
hyperplane section of $(n-1)$-dimensional volume at least $c$. Here and
below, a convex body is a compact convex set with nonempty interior.
An equivalent formulation is that the isotropic constant
\begin{equation}\label{eq:isotropic-constant}
L_K=\frac{\det(\cov(K))^{1/(2n)}}{|K|^{1/n}}
\end{equation}
is bounded by a universal constant, independently of $n$ and $K$.
Here $\cov(K)$ denotes the covariance of the uniform probability measure on $K$,
and $|K|$ is the volume of $K$. The quantity $L_K$ is affine invariant.
When the uniform probability measure on $K$ is isotropic, i.e., when it has barycenter zero
and covariance matrix equal to the identity, we have $L_K = |K|^{-1/n}$.
See Milman and Pajor~\cite{MP} for background on the isotropic constant
and its relation to hyperplane sections.

The slicing theorem strengthens Milman's theory of $M$-ellipsoids~\cite{Milman}.
One formulation of that theory asserts the existence of an ellipsoid
$\mathcal E$ with $|\mathcal E|=|K|$ such that
$N(K,\mathcal E)\leq C^n$, where $N(K,\mathcal E)$ is the minimum
number of translates of $\mathcal E$ needed to cover $K$ and where $C$ is
a universal constant. We say that $K$ is in $M$-position when a
Euclidean ball has this property.
Elementary volumetric considerations show that, when $K$ is isotropic
and $D$ is the centered Euclidean ball with $|D|=|K|$,
\[
N(K,D)\leq\exp(C'nL_K),
\]
where $C'$ is universal; see, e.g., \cite[Lemma 3.4]{GPT}.
A uniform upper bound on $L_K$ therefore makes $D$ an $M$-ellipsoid
of $K$, so isotropic position is an $M$-position.
Thus a proof of slicing which does not appeal to $M$-ellipsoids also
recovers their existence.

The slicing problem originated in Bourgain's work in the 1980s~\cite{Bourgain}.
His bound of order $n^{1/4}\log n$~\cite{Bourgain91} was improved to $n^{1/4}$ by
Klartag~\cite{K06} using exponential tilts. Subsequent progress was closely connected to two
other problems concerning log-concave measures. The
Kannan--Lov\'asz--Simonovits (KLS) conjecture~\cite{KLS}, in its
Poincar\'e formulation, asserts that every isotropic log-concave
probability measure $\mu$ satisfies, for every locally Lipschitz function $f$,
\begin{equation}\label{eq:KLS}
\var_\mu(f)\leq C\int_{\R^n}|\nabla f(x)|^2\,d\mu(x)
\end{equation}
with a universal $C$. Ball and Nguyen~\cite{BN} showed that this
implies slicing. The thin-shell conjecture asks for a universal bound
on $\E(|X|-\sqrt n)^2$ when $X$ is isotropic and log-concave.
Equivalently, up to universal constants, it asks for
$\var(|X|^2)\leq Cn$; see~\cite{GM}. It is a consequence of
\eqref{eq:KLS}, applied to the Euclidean norm (or its square), and
Eldan and Klartag~\cite{EK} proved that it too implies slicing.

Eldan's stochastic localization~\cite{Eldan} provided a powerful tool
for these questions. Chen~\cite{Chen} obtained a subpolynomial bound
for KLS, and hence for slicing. Klartag and Lehec~\cite{KLpoly}
subsequently obtained polylogarithmic bounds. Building on a covariance estimate
of Guan~\cite{Guan}, Klartag and Lehec~\cite{KLslicing} proved the
slicing theorem. Bizeul~\cite{Bizeul} subsequently gave a proof through
small-ball estimates. The first proof uses $M$-ellipsoids directly;
in the second, they enter through the reduction of Dafnis and
Paouris~\cite{DP} from small-ball estimates to slicing.

Klartag and Lehec~\cite{KLthin} then proved the thin-shell conjecture,
which gives another route to slicing through~\cite{EK}. The latter
implication uses the reverse Santal\'o inequality. More recently,
Chen and Klartag~\cite{CK} presented an AI-assisted proof of the sharp
bound $\var(|X|^2)\leq 8n$. Their method uses optimal transport and
log-concave moment measures rather than stochastic localization.
As observed by Bizeul~\cite{BizeulMM}, it yields a similar estimate
for all quadratic forms, not just the squared Euclidean norm.
We return to this in Section~\ref{sec:preliminaries}.
Brazitikos~\cite{Brazitikos} developed a deterministic approach based on
quadratic tilts, and Brazitikos and Milman~\cite{BM} used central
Gaussian tilts to pass from thin-shell bounds to small-ball estimates.
The geometric passage from these estimates to slicing uses either
the Dafnis--Paouris reduction or the alternative completion in
\cite{Brazitikos}, which invokes reverse Santal\'o. Thus these approaches
retain an input from $M$-ellipsoid theory or the closely related reverse
Santal\'o theory. The KLS conjecture itself was recently resolved in
three contemporaneous works by Balasubramanian and Kasiviswanathan~\cite{BaKa},
Bizeul, Klartag and Lehec~\cite{BKL}, and Song and Zhang~\cite{SZ}.

The purpose of this note is to give a short passage from
\eqref{eq:strong-thin-shell} to slicing which uses neither of these
geometric inputs and gives an explicit constant.

\begin{theorem}[Slicing]\label{thm:slicing}
For every convex body $K\subset\R^n$,
\[
L_K\leq \sqrt{\frac5\pi}.
\]
\end{theorem}

The argument can be summarized as follows.
The strong thin-shell bound controls the third-moment
tensor of every isotropic log-concave measure. Along stochastic
localization, this gives a lower bound on the expected logarithmic
determinant of the conditional covariance. A change of variables then
relates this expected log-determinant to the isotropic constant.

Section~\ref{sec:preliminaries} recalls the cumulant bound and the
stochastic localization identities used in the proof.
Section~\ref{sec:slicing} proves Theorem~\ref{thm:slicing}.

Finally, an OpenAI preprint~\cite{OpenAI} claims the sharp inequality
$L_K\leq L_{\Delta_n}$, where $\Delta_n$ is an $n$-dimensional simplex.
We do not assess that claimed proof here and use none of its results.
Observe, however, that since
\[
L_{\Delta_n} = \frac{n!^{1/n}}{ (n+1)^{\frac{n+1}{2n}} \sqrt{n+2} } \to \frac 1\e
\]
as $n$ tends to infinity, our constant $\sqrt{5/\pi}$ exceeds the
asymptotically optimal value by a factor $\e\sqrt{5/\pi}\simeq3.43$.

\section{Preliminaries}
\label{sec:preliminaries}

\paragraph{Third-order cumulant bound.}
Chen and Klartag~\cite{CK} proved that every isotropic log-concave
random vector $X$ in $\R^n$ satisfies
\[
\var(|X|^2)\leq 8n.
\]
Equality holds when the coordinates of $X$ are independent centered
exponential random variables of variance one. Their proof is based on
optimal transport, more specifically on log-concave moment measures.
The following anisotropic form was proved by Letwin~\cite{Letwin}
with an intricate argument. It also follows from a simple modification
of Chen and Klartag's argument; see Bizeul~\cite{BizeulMM}.
For every isotropic log-concave random vector $X$ in $\R^n$ and every
symmetric matrix $A$,
\begin{equation}\label{eq:strong-thin-shell}
\var( \langle AX,X\rangle)\leq 8\tr\bigl[A^2\bigr].
\end{equation}

Following~\cite{BizeulMM}, we obtain a universal bound on the third-order
cumulants of isotropic log-concave measures on $\R^n$. Given such a measure $\mu$, let
\[
\kappa_3 (\mu) = \sup_{\theta\in \SS^{n-1}} \left\{  \Vert \E [ \langle X,\theta\rangle  X \otimes X ] \Vert \right\} ,
\]
where $X \sim \mu$ and the norm is the Euclidean (Hilbert--Schmidt)
norm of the matrix $\E [\langle X,\theta\rangle X\otimes X]$.
\begin{theorem} \label{main2}
For every positive integer $n$ and every isotropic
log-concave probability measure $\mu$ on $\R^n$, we have
\[
\kappa_3 (\mu) \leq 2\sqrt 2.
\]
\end{theorem}
\begin{proof}
Fix $\theta\in\SS^{n-1}$ and let $X\sim\mu$. Set
$H=\E[\langle X,\theta\rangle X\otimes X]$, a symmetric matrix.
By isotropy, $\langle X,\theta\rangle$ has expectation zero and
variance one. Cauchy--Schwarz and~\eqref{eq:strong-thin-shell} give
\[
\tr H^2
= \E [\langle X,\theta\rangle \langle H X , X \rangle]
\leq \left( \var(\langle HX,X\rangle) \right)^{1/2}
\leq  ( 8 \cdot \tr H^2 )^{1/2} ,
\]
hence $( \tr H^2 )^{1/2} \leq 2 \sqrt 2$. Taking the supremum over
$\theta$ proves the result.
\end{proof}
The parameter
\[
\kappa_3^*(n):=\sup\{\kappa_3(\mu):\mu\text{ is isotropic and
log-concave on }\R^n\}
\]
was introduced by Eldan~\cite{Eldan} in connection with stochastic
localization and KLS. With $C_P(\mu)$ denoting the best constant in
$\var_\mu(f)\leq C_P(\mu)\int|\nabla f(x)|^2\,d\mu(x)$, the estimate
of Klartag~\cite{Klog} gives
\[
C_P(\mu)\leq C\sqrt{\log n}\,\kappa_3^*(n)
\]
for isotropic log-concave $\mu$. Thus Theorem~\ref{main2} recovers
the $O((\log n)^{1/4})$ bound for the KLS constant;
see~\cite{Letwin}. This consequence is superseded by the
recent dimension-free results in~\cite{BaKa,BKL,SZ}.

In fact, an extension of Theorem~\ref{main2} to higher-order cumulants
is common to these proofs of the KLS conjecture. For $m\geq2$, define
\[
K_m(\mu)=\left.\nabla^m\log\int_{\R^n}\e^{\langle z,x\rangle}\,d\mu(x)
\right|_{z=0}
\]
and
\[
\kappa_m(\mu)=\sup_{\theta\in\SS^{n-1}}
\left(\sum_{i_1,\ldots,i_{m-1}=1}^n
K_m(\mu)(e_{i_1},\ldots,e_{i_{m-1}},\theta)^2\right)^{1/2},
\]
where $(e_i)$ is an orthonormal basis. For isotropic $\mu$, we have
$\kappa_2(\mu)=1$, and the definition of $\kappa_3$ agrees with the one
above. The three works~\cite{BaKa,BKL,SZ} establish, in different
forms, the bound
\[
\kappa_m(\mu)\leq (Cm)^m,\qquad m\geq2,
\]
with a universal constant $C$, for every isotropic log-concave measure
$\mu$. These higher-order estimates are used to obtain a universal
bound on $C_P(\mu)$.

\paragraph{Stochastic localization.}

We briefly recall the stochastic localization construction and refer to~\cite[Sections 6 and 7.1--7.2]{KLsurvey} for details.

Let $K \subset \R^n$ be a convex body.
Let $X$ be a random vector uniformly distributed on $K$,
and let $(W_t)_{t \geq 0}$ be a standard Brownian motion in $\R^n$,
independent of $X$. Set
\[
\theta_t = tX + W_t ,
\]
and let $\mathcal F_t=\sigma(\theta_s:0\leq s\leq t)$ be the observation
filtration. Let $\mu_t$ be the conditional law of $X$ given $\mathcal F_t$,
and let $a_t$ be the barycenter of $\mu_t$:
\[
a_t = \int_{\R^n} x\,d\mu_t(x) = \E [ X \mid  \mathcal F_t ].
\]
A direct computation shows that $\mu_t$ has density
\[
p_t(x)=\frac{d\mu_t}{d\mu_0}(x)
=\exp\left(\langle\theta_t,x\rangle-\frac t2|x|^2
-\Lambda_t(\theta_t)\right),
\]
with respect to the uniform probability measure $\mu_0$ on $K$, where
\[
\Lambda_t (\theta) =
\log \left[ \E \exp \left( \langle \theta, X \rangle - \frac t2 \vert X\vert^2 \right) \right]
\]
is the log-Laplace transform of the Gaussian tilt.
Differentiating under the integral sign gives
\[
a_t = \nabla\Lambda_t(\theta_t), \qquad \nabla^2 \Lambda_t ( \theta_t ) = \cov(\mu_t) =: A_t .
\]
Filtering theory shows that the process $(B_t)$ given by
\[
B_t=\theta_t-\int_0^t a_s\,ds , \qquad t > 0,
\]
is a standard Brownian motion with respect to $(\mathcal F_t)$; see, e.g., \cite{KLthin}. Thus $(\theta_t)$ satisfies the stochastic differential
equation
\[
d\theta_t=dB_t+a_t\,dt = dB_t + \nabla \Lambda_t ( \theta_t ) \,dt .
\]
For any bounded test function $\phi$, the process
\[
M_t:=\int_{\R^n} \phi(x)\,d\mu_t(x)=\E[\phi(X)\mid\mathcal F_t]
\]
is a martingale and admits the It\^o representation
\[
M_t-M_0=\int_0^t\langle v_s,dB_s\rangle,
\qquad v_t=\int_{\R^n} \phi(x)(x-a_t)\,d\mu_t(x).
\]
Applying this identity to linear and quadratic functions and using
It\^o's formula gives the covariance equation
\[
dA_t=\sum_{i=1}^n H_{i,t}\, \langle dB_t,e_i\rangle -A_t^2\,dt,
\]
where $e_1,\dotsc,e_n$ is the canonical basis of $\R^n$ and
\begin{equation}\label{eq_Hi}
H_{i,t}=\int_{\R^n} (x-a_t)\otimes(x-a_t) \langle x-a_t , e_i \rangle\,d\mu_t(x), \qquad i = 1 ,\dotsc,n.
\end{equation}
It\^o calculus then gives, for any smooth function $f\colon (0,\infty) \to \R$,
\begin{equation}\label{eq:matrix-ito}
\begin{split}
d\tr f(A_t)
={}&\sum_{i=1}^n\tr\bigl[f'(A_t)H_{i,t}\bigr]\,\langle dB_t,e_i\rangle
-\tr\bigl[f'(A_t)A_t^2\bigr]\,dt\\
&+\frac12\sum_{i,j,k=1}^n
\frac{f'(\lambda_{i,t})-f'(\lambda_{j,t})}
{\lambda_{i,t}-\lambda_{j,t}}\,\xi_{ijk,t}^2\,dt,
\end{split}
\end{equation}
where the $\lambda_{i,t}$ are the eigenvalues of $A_t$ and
$\xi_{ijk,t}$ are the coordinates of the tensor
\[
\int_{\R^n}  (x-a_t)^{\otimes 3}\,d\mu_t(x)
\]
in the corresponding orthonormal eigenbasis.
The divided difference is
interpreted as $f^{\prime\prime}(\lambda_{i,t})$ when the two eigenvalues agree.

Under suitable integrability conditions, the stochastic integral is a
genuine martingale. For instance, when $f(x)=x$, this follows from the
boundedness of $K$, since all centered moments of $\mu_t$ are then
uniformly bounded. Taking expectations gives
\[
\frac d{dt}\E\tr A_t=-\E\tr A_t^2,
\]
so $\E\tr A_t$ is decreasing. Since $A_0=\cov(K)$, in the isotropic case we obtain
\begin{equation}\label{eq_trA_t}
\E \tr A_t \leq n, \qquad t \geq 0 .
\end{equation}

In the next section we apply~\eqref{eq:matrix-ito} to $f(x)=\log x$.
We also justify the martingale property in this case, despite the
singularity of $f'$ at zero.

\section{Proof of slicing}
\label{sec:slicing}

The quantity $L_K$ is affine invariant, so we can assume
that $K$ is isotropic, in the sense that the barycenter
is at $0$ and the covariance matrix of the uniform probability measure on $K$ is the identity.

Write $S(Y)=-\int_{\R^n}p(x)\log p(x)\,dx$ for the differential entropy
of a random vector $Y$ with density $p$. For $t>0$, the random vector
$a_t$ is supported on $K$, so its entropy satisfies
\[
S ( a_t ) \leq S ( X ) = \log \vert K \vert .
\]
Recall that $a_t=\nabla\Lambda_t(\theta_t)$ and
$\nabla^2\Lambda_t(\theta_t)=A_t$.
The Hessian of $\Lambda_t$ is positive definite, so $\nabla\Lambda_t$ is
one-to-one and is a diffeomorphism from $\R^n$ to the interior of $K$.  Changing variables gives
\[
S (a_t) = S ( \theta_t ) + \E \log \det A_t.
\]
The fact that conditioning decreases entropy gives
$S(tX+W_t)\geq S(tX+W_t\mid X)=S(W_t)$, hence
\[
S ( \theta_t ) = S ( t X + W_t ) \geq S ( W_t ) = \frac n2 \log ( 2 \pi \e) + \frac n2 \log t .
\]
Putting everything together, we get
\begin{equation}
\label{eq_zou}
\log \vert K \vert \geq \frac n2 \log ( 2\pi \e ) + \frac n2 \log t +
\E \log \det A_t.
\end{equation}
We next estimate $\log\det A_t$. Note that
\[
\log \det A_t = \tr \log A_t = \sum_{k} \log \lambda_{k, t} ,
\]
where $0 < \lambda_{1,t} \leq \ldots \leq \lambda_{n,t}$ are the eigenvalues of $A_t$,
repeated according to their multiplicity. By isotropy, $A_0=\cov(K)=\id$.
For $f(x)=\log x$, formula~\eqref{eq:matrix-ito} gives
\[
\log\det A_t=M_t+\int_0^t\mathcal D_s\,ds,
\]
where $M_t$ is a continuous local martingale with $M_0=0$, and
\[
\begin{split}
\mathcal D_t
&=-\tr A_t-\frac12\sum_{ijk}
\frac{\xi_{ijk,t}^2}{\lambda_{i,t}\lambda_{j,t}}\\
&=-\tr A_t-\frac12\sum_{ijk}\lambda_{k,t}\chi_{ijk,t}^2.
\end{split}
\]
Here
$\chi_{ijk,t}=\xi_{ijk,t}/\sqrt{\lambda_{i,t}\lambda_{j,t}\lambda_{k,t}}$
are the coordinates of the third-moment tensor of the centered,
isotropic image of $\mu_t$, in the basis for which $A_t$ is diagonal.
Since this image is log-concave,
Theorem~\ref{main2} gives
\[
\sum_{ij}\chi_{ijk,t}^2\leq8
\]
for every $k$ and $t$, almost surely. Consequently,
\[
\mathcal D_t \geq - \tr A_t - 4\sum_k \lambda_{k,t} = - 5\tr A_t ,
\]
and hence
\begin{equation} \label{eq_logdet}
\log\det A_t\geq M_t-5\int_0^t\tr A_s\,ds.
\end{equation}
Moreover, $M$ is a genuine martingale with mean zero. Indeed, by~\eqref{eq:matrix-ito},
\[
dM_t=\langle v_t,dB_t\rangle,
\qquad v_t=\sum_{i=1}^n\tr(A_t^{-1}H_{i,t})e_i,
\]
where $(e_i)$ is the canonical basis and the matrices $H_{i,t}$ are
defined in~\eqref{eq_Hi}. Computing in an orthonormal eigenbasis of
$A_t$ gives
\[
|v_t|^2=\sum_k\lambda_{k,t}\left(\sum_i\chi_{iik,t}\right)^2.
\]
By Cauchy--Schwarz and the bound $\sum_{ij}\chi_{ijk,t}^2\leq8$,
\[
\left(\sum_i\chi_{iik,t}\right)^2
\leq n\sum_i\chi_{iik,t}^2
\leq n\sum_{ij}\chi_{ijk,t}^2
\leq 8n.
\]
Thus $|v_t|^2\leq 8n\tr A_t$, and by~\eqref{eq_trA_t},
\[
\E\int_0^t |v_s|^2\,ds
\leq 8n\int_0^t \E\tr A_s\,ds
\leq 8n^2t < +\infty.
\]
Therefore $M$ is a square-integrable martingale on every finite time interval. In particular, $\E M_t=0$ for every $t\geq0$.
Taking expectations in~\eqref{eq_logdet}
and using~\eqref{eq_trA_t}, we obtain
\[
\E\log\det A_t\geq-5\int_0^t\E\tr A_s\,ds\geq-5nt.
\]
Substituting into~\eqref{eq_zou} yields
\[
\log \vert K \vert \geq \frac n2 \log ( 2\pi \e ) + \frac n2 \log t - 5nt ,
\]
for every positive $t$. Choosing the maximizing value $t=1/10$, we get
\[
L_K = \frac 1{ |K|^{1/n} } \leq \sqrt{\frac5\pi} .
\]

\begin{remark}
The same argument yields a better constant.
First, writing $m(t)=\E\tr A_t$, the covariance equation gives
\[
m'(t)=-\E\tr A_t^2\leq-\frac{m(t)^2}{n},
\qquad  m(t)\leq\frac{n}{1+t}.
\]
One can also retain the exact entropy contribution through the
I--MMSE identity of Guo, Shamai and Verd\'u~\cite{GSV}:
\[
S(\theta_t)=\frac n2\log(2\pi\e t)+\frac12\int_0^t m(s)\,ds.
\]
Combining these refinements yields
\[
L_K\leq\left(\frac{(1+t)^9}{2\pi\e t}\right)^{1/2},
\]
and choosing $t=1/8$ gives
$L_K\leq\bigl(9^9/(2\pi\e\,8^8)\bigr)^{1/2}\simeq1.163$.
These improvements remain far from the simplex value.
\end{remark}

\paragraph{Use of AI tools.}
All mathematical content is due to the authors, with the
exception of the very last remark, involving the I--MMSE identity,
which was found by an AI agent. AI tools were also used to
draft parts of the exposition from outlines and detailed instructions
provided by the authors, and to refine the wording.

\bigskip
\noindent School of Mathematical Sciences, Tel Aviv University, Tel Aviv, Israel; and \\
Department of Mathematics, Weizmann Institute of Science, Rehovot, Israel. \\
\textit{e-mail:} \texttt{klartagb@tau.ac.il}

\bigskip
\noindent Institut de Mathématiques de Jussieu--Paris Rive Gauche,
Équipe d'Analyse Fonctionnelle, Sorbonne Université, Paris, France. \\
\textit{e-mail:} \texttt{joseph.lehec@imj-prg.fr}

\end{document}